\documentclass[12pt,reqno]{amsart}
\usepackage{enumerate, latexsym, amsmath, amsfonts, amssymb, amsthm, graphicx, color}
\def\pmod #1{\ ({\rm{mod}}\ #1)}
\def\Z{\Bbb Z}
\def\N{\Bbb N}

\def\l{\left}
\def\r{\right}
\def\bg{\bigg}
\def\({\bg(}
\def\){\bg)}
\def\t{\text}
\def\f{\frac}

\def\ls{\leqslant}
\def\gs{\geqslant}

\def\eq{\equiv}

\theoremstyle{plain}
\newtheorem{theorem}{Theorem}

\newtheorem{lemma}{Lemma}

\newtheorem{conjecture}{Conjecture}
\theoremstyle{definition}

\theoremstyle{remark}

\makeatletter
\@namedef{subjclassname@2010}{%
  \textup{2010} Mathematics Subject Classification}
\makeatother
 \vspace{4mm}
\begin{document}

\hbox{Preprint}
\medskip

\title{On Practical Numbers of Some Special Forms}

\author{Li-Yuan Wang}
\address{Department of Mathematics, Nanjing University\\
Nanjing, 210093, People's Republic of China}
\email{dg1721010@smail.nju.edu.cn}

\author{Zhi-Wei Sun }
\address{Department of Mathematics, Nanjing University\\
Nanjing, 210093, People's Republic of China}
\email{zwsun@nju.edu.cn}

\begin{abstract} In this paper we study practical numbers of some special forms.
For any integers $b\gs0$ and $c>0$, we show that if $n^2+bn+c$ is practical for some integer $n>1$, then
there are infinitely many nonnegative integers $n$ with $n^2+bn+c$ practical. We also prove that there are infinitely many practical numbers of the form $q^4+2$ with $q$ practical, and that
there are infinitely many practical Pythagorean triples $(a,b,c)$ with $\gcd(a,b,c)=6$ (or $\gcd(a,b,c)=4$).
\end{abstract}

\thanks{
\newline \indent The second author is the corresponding author, and supported by the National Natural Science
Foundation of China (grant 11571162).}

\maketitle
\section{Introduction}	
\setcounter{lemma}{0}
\setcounter{theorem}{0}
\setcounter{corollary}{0}
\setcounter{remark}{0}
\setcounter{equation}{0}
\setcounter{conjecture}{0}

A positive integer $m$ is called a {\it practical number} if each $n=1,\ldots,m$ can be written as the sum of some distinct divisors of $n$. This concept was introduced by Srinivasan \cite{name-of-practical} who noted that any practical number greater than $2$ must be divisible by $4$ or $6$. In 1954, Stewart \cite{structure-theorem} obtained the following structure theorem for practical numbers.
\begin{theorem} Let $p_1<\ldots<p_k$ be distinct primes and let $a_1,\ldots,a_k\in\Z^+=\{1,2,3,\ldots\}$. Then
$m=p_{1}^{a_{1}}p_{2}^{a_{2}} \cdots p_{k}^{a_{k}}$ is practical if and only if $p_{1}=2$ and
$$p_{j}-1\ls \sigma (p_{1}^{a_{1}}p_{2}^{a_{2}} \cdots p_{j-1}^{a_{j-1}})\ \ \t{for all}\ 1<j\ls k,$$
where $\sigma(n)$ denotes the sum of the positive divisors of $n$.
\label{structure}
\end{theorem}

It is interesting to compare practical numbers with primes. All practical numbers are even except 1 while all primes are odd except 2. Moreover, if $P(x)$ denotes the number of practical numbers not exceeding $x$, then there is a positive constant $c$ such that
\begin{equation}
 P(x)\sim \frac{cx}{\log x}\quad\t{as}\ x\to\infty,
 \label{zero}
 \end{equation}
which was established in \cite{distribution}. This is quite similar to the Prime Number Theorem.

 Inspired by the famous Goldbach conjecture and the twin prime conjecture, Margenstern \cite{french} conjectured that every positive even integer is the sum of two practical numbers and that there are infinitely many practical numbers $m$ with $m-2$ and $m+2$ also practical. Both conjectures were confirmed by Melfi \cite{two-conjecture} in 1996. Guo and Weingartner \cite{GW} proved in 2018 that for any odd integer $a$ there are infinitely many primes $p$ with $p+a$ practical.
 An open conjecture of Sun \cite[Conjecture 3.38] {S17} states that any odd integer greater than one can be written as the sum of a prime and a practical number.

 Whether there are infinitely many primes of the form $x^2+1$ with $x\in\Z$ is a famous unsolved problem in number theory. Motivated by this, in 2017 Sun \cite[A294225]{S} conjectured that there are infinitely many positive integers $q$
  such that $q$, $q+2$ and $q^2+2$ are all practical, which looks quite challenging. Thus, it is natural to study for what $a,b,c\in\Z^+$ there are infinitely many
  practical numbers of the form $an^2+bn+c$. Note that if $a\eq b\pmod 2$ and $2\nmid c$ then $an^2+bn+c$ is odd for any $n\in\N=\{0,1,2,\ldots\}$ and hence
  $an^2+bn+c$ cannot take practical values for infinitely many $n\in\N$.

Based on our computation we formulate the following conjecture.

\begin{conjecture} Let $a,b,c$ be positive integers with $2\nmid ab$ and $2\mid c$. Then there are infinitely many
$n\in\N$ with $an^2+bn+c$ practical. Moreover, in the case $a=1$, there is an integer $n$ with $1<n\ls\max\{b,c\}$ such that
$n^2+bn+c$ is practical.
\end{conjecture}

Though we are unable to prove this conjecture fully, we make the following progress.

\begin{theorem}\label{Th1.1}
Let $b\in\N$ and $c\in\Z^+$. If $n^2+bn+c$ is practical for some integer $n>1$, then
there are infinitely many $n\in\N$ with $n^2+bn+c$ practical.
\label{quadratic}
\end{theorem}

If $1\ls b\ls 100$ and $1\ls c\ls100$ with $2\nmid b$ and $2\mid c$, then we can easily find $1<n\ls \max\{b,c\}$
with $n^2+bn+c$ practical. For example,  $n^2+n+2$ with $n=2$ is practical. For each positive even number $b\ls 20$ we make the set
$$S_b:=\{1\ls c\ls 100:\ n^2+bn+c\ \t{is practical for some}\ n=2,\ldots,20000\}$$
explicit:
\begin{align*}S_0=&\{1\ls c\ls 100:\ c\not\eq 1,10\pmod{12}\ \t{and}\ c\not=43,67,93\},
\\S_2=&\{1\ls c\ls 100:\ c\not\eq2,11\pmod{12}\ \t{and}\ c\not=44,68,94\},
\\S_4=&\{1\ls c\ls 100:\ c\not\eq 2,5\pmod{12}\ \t{and}\ c\not=47,71,97\},
\\S_6=&\{1\ls c\ls 100:\ c\not\eq7,10\pmod{12}\ \t{and}\ c\not=52,76\},
\\S_8=&\{1\ls c\ls 100:\ c\not\eq2,5\pmod{12}\ \t{and}\ c\not=59,83\},
\\S_{10}=&\{1\ls c\ls 100:\ c\not\eq 2,11\pmod{12}\ \t{and}\ c\not=68,92\},
\\S_{12}=&\{1\ls c\ls 100:\ c\not\eq1,10\pmod{12}\ \t{and}\ c\not=79\},
\\S_{14}=&\{1\ls c\ls 100:\ c\not\eq2,11\pmod{12}\ \t{and}\ c\not=92\},
\\S_{16}=&\{1\ls c\ls 100:\ c\not\eq 2,5\pmod{12}\},
\\S_{18}=&\{1\ls c\ls 100:\ c\not\eq7,10\pmod{12}\},
\\S_{20}=&\{1\ls c\ls 100:\ c\not\eq2,5\pmod{12}\}.
\end{align*}
For example, applying Theorem \ref{Th1.1} with $b=20$, we see that for any $c=1,\ldots,100$ with $c\not\eq2,5\pmod{12}$
there are infinitely many $n\in\N$ with $n^2+20n+c$ practical. It is easy to see that if $c$ is congruent to $2$ or $5$ modulo $12$ then $n^2+20n+c$ is not practical for any integer $n\gs2$.

By Theorem \ref{Th1.1} and the fact $2\in S_0$, there are infinitely many $n\in\N$ with $n^2+2$ practical. Moreover, we have the following stronger result.

\begin{theorem}\label{Th1.2} $2^{35\times3^k+1}+2$ is practical for every $k=0,1,2,\ldots$. Hence
there are infinitely many practical numbers $q$ with $q^4+2$ also practical.
\label{infinity}
\end{theorem}

We prove Theorem \ref{Th1.2} by modifying Melfi's cyclotomic method in \cite{two-conjecture}.

We now turn to Pythagorean triples involving practical numbers, and call
a Pythagorean triple $(a,b,c)$ with $a,b,c$ all practical a {\it practical Pythagorean triple}.
Obviously, there are infinitely many practical Pythagorean triples. In fact, if $a^2+b^2=c^2$ with $a,b,c$ positive integers then $({2^ka})^2+({2^kb})^2=({2^kc})^2$ for all $k=0,1,2,\ldots$. By Theorem \ref{structure}, $2^ka$, $2^kb$ and  $2^kc$ are all practical if $k$ is large enough.

 Our following theorem was originally conjectured by Sun \cite[A294112]{S}.

\begin{theorem}\label{Th1.3} Let $d$ be $4$ or $6$. Then
there are infinitely many practical Pythagorean triples $(a,b,c)$ with $\gcd(a,b,c)=d$.
\label{Pythagorean}
\end{theorem}

We are going to show Theorems \ref{Th1.1}-\ref{Th1.3} in the next section.
\section{Proofs of Theorems 1.2-1.4}

\begin{lemma} \label{Lem2.1} Let $m$ be any practical number. Then $mn$ is practical for every $n=1,\ldots,\sigma(m)+1$.
 In particular, $mn$ is practical for every $1\ls n\ls 2m$.
\label{twotimes}
\end{lemma}
This lemma follows easily from Theorem \ref{structure}; see \cite{two-conjecture} for details. Note that if $m>1$ is practical then
$m-1$ can be written as the sum of some divisors of $m$ and hence $(m-1)+m\ls\sigma(m)$.
\medskip

\noindent{\it Proof of Theorem \ref{Th1.1}}. Set $f(n)=n^2+bn+c$. It is easy to verify that
$$f(n+f(n))=f(n)(f(n)+2n+b+1).$$
Note that
$$f(n)-(2n+b+1)=n(n-2)+b(n-1)+c-1\gs0.$$
If $n\gs2$ is an integer with $f(n)$ practical, then $f(n+f(n))=f(n)(f(n)+2n+b+1)$ is also practical by Lemma \ref{Lem2.1} and the inequality
$$f(n)+2n+b+1\ls 2f(n).$$ So the desired result follows.  \qed
\medskip

For a positive integer $m$, the cyclotomic polynomial $\Phi_m(x)$ is defined by
$$\Phi_m(x):=\prod_{a=1\atop\gcd(a,m)=1}^m\l(x-e^{2\pi ia/m}\r).$$
Clearly,
\begin{equation}\label{cyc}x^n-1=\prod_{d\mid n}\Phi_{d}(x)\quad\t{for all}\ n=1,2,3,\ldots.
\end{equation}
\medskip

\noindent{\it Proof of Theorem \ref{infinity}}.
Write $m_k=2^{35\times3^k+1}+2$ for $k=0,1,2,\ldots$. Note that
$m_{2k}=q_k^4+2$ with $q_k=2^{(35\times9^k+1)/4}$ practical. So it suffices to prove that $m_k$
is practical for every $k=0,1,2,\ldots$.

Via a computer we find that
$$m_{0}=2^{36}+2,\ m_{1}=2^{106}+2,\ m_{2}=2^{316}+2$$
 are all practical.

Now assume that $m_{k}$ is practical for a fixed integer $k\gs 2$.  For convenience, we write $x$ for $2^{3^k}$. Then
$$x\gs 2^9=512,\ m_k=2(x^{35}+1)\ \t{and}\ m_{k+1}=2(x^{105}+1).$$
 In view of \eqref{cyc},
\begin{equation}\label{105}
\f{x^{210}-1}{x^{105}-1}=\f{x^{70}-1}{x^{35}-1}\Phi_{6}(x)\Phi_{30}(x)\Phi_{42}(x)\Phi_{210}(x).
\end{equation}
Since $x\gs512$, we have
\begin{equation}
 \frac{x^2}{2}< \Phi_{6}(x)=x^2-x+1<x^2.
 \label{phi6}
\end{equation}
Clearly,
$$x^7>x^3\f{x^3-1}{x-1}=x^5+x^4+x^3$$
and
$$x^8>2x^7\gs x^7+x+1.$$
Thus
\begin{equation}
 x^8< \Phi_{30}(x)=x^8+x^7-x^5-x^4-x^3+x+1<2x^8
 \label{phi30}
\end{equation}
Similarly, for
$$\Phi_{42}(x) =x^{12}+x^{11}-x^9-x^8+x^6-x^4-x^3+x+1$$
and
\begin{align*}
\Phi_{210}(x)=&x^{48}-x^{47}+x^{46}+x^{43}-x^{42}+2 x^{41}-x^{40}+x^{39}+x^{36}
\\&-x^{35}+x^{34}-x^{33}+x^{32}-x^{31}-x^{28}-x^{26}-x^{24}-x^{22}
\\&-x^{20}-x^{17}+x^{16}-x^{15}+x^{14}-x^{13}+x^{12}+x^9-x^8
\\&+2x^7-x^6+x^5+x^2-x+1,
\end{align*}
we can prove that
\begin{equation}\label{phi42}
 x^{12}< \Phi_{42}(x)<2x^{12}\ \t{and}\
 \Phi_{210}(x)<x^{48}.
\end{equation}
Combining (\ref{phi6}), (\ref{phi30}) and (\ref{phi42}), we get
\begin{equation}
\frac{x^{22}}{2}<\Phi_{6}(x)\Phi_{30}(x)\Phi_{42}(x)<4x^{22}
\label{inequality}
\end{equation}
and hence $\Phi_{6}(x)\Phi_{30}(x)\Phi_{42}(x)<4(x^{35}+1)$. Thus, by Lemma \ref{twotimes} and the induction hypothesis we obtain that
$$2(x^{35}+1)\Phi_{6}(x)\Phi_{30}(x)\Phi_{42}(x)$$
is practical.

By (\ref{inequality}),
$$2(x^{35}+1)\Phi_{6}(x)\Phi_{30}(x)\Phi_{42}(x)>x^{57}>x^{48}.$$
So, applying (\ref{phi42}) and Lemma \ref{twotimes}, we conclude that
$$2(x^{35}+1)\Phi_{6}(x)\Phi_{30}(x)\Phi_{42}(x)\Phi_{210}(x)$$
is practical. In view of \eqref{105}, this indicates that $m_{k+1}$ is practical. This completes the proof. \qed

\begin{lemma} {\rm (Melfi \cite{two-conjecture})}
For every $k \in \N$, both $ 2(3^{3^k\cdot 70}-1) $ and $ 2(3^{3^k\cdot 70}+1)$ are practical numbers.
\label{bothpractical}
\end{lemma}

\medskip
\noindent{\it Proof of Theorem} 1.4. (i) We first consider the case $d=4$. For each $k=0,1,2,\ldots$, define
$$a_k=2(3^{3^k\cdot 70}-1),\ b_k=4\cdot 3^{3^k\cdot 35},\ \t{and}\ c_k=2(3^{3^k\cdot 70}+1).$$
It is easy to see that $a_k^2+b_k^2=c_k^2$ and $\gcd(a_k,b_k,c_k)=4$.
By Lemma \ref{bothpractical}, $a_k$ and $c_k$ are both practical. Theorem \ref{twotimes} implies that $b_k$ is practical. This proves Theorem 1.4 for $d=4$.

(ii) Now we handle the case $d=6$. For any $k=0,1,2,\ldots$, define
$$x_k=3(3^{3^k\cdot 70}-1),\ y_k=6\cdot 3^{3^k\cdot 35},\ \t{and}\ z_k=3(3^{3^k\cdot 70}+1).$$
 Then $x_k^2+y_k^2=z_k^2$ and $\gcd(x_k,y_k,z_k)=6$. Note that
 $y_k$ is practical for any $k=0,1,2,\ldots$ by Theorem \ref{twotimes}.

 Now it remains to show by induction that $x_k$ and $z_k$ are practical for all $k=0,1,2,\ldots$.
 Via a computer, we see that $x_0=3^{71}-3$ and $z_0=3^{71}+3$ are practical numbers.
 Suppose that $x_k$ and $z_k$ are practical for some nonnegative integer $k$. Then
 \begin{equation}\label{2.8}
 x_{k+1}=3(3^{3^{k+1}\cdot 70}-1)=x_k(3^{3^k\cdot 70}-3^{3^k\cdot 35}+1)(3^{3^k\cdot 70}+3^{3^k\cdot 35}+1)
 \end{equation}
 and
 \begin{equation}\label{2.9}
z_{k+1}=3(3^{3^{k+1}\cdot 70}+1)= z_k\Phi_{12}(3^{3^k})\Phi_{60}(3^{3^k})\Phi_{84}(3^{3^k})\Phi_{420}(3^{3^k}).
 \end{equation}
In view of (\ref{2.8}), by applying Lemma \ref{twotimes} twice, we see that $x_{k+1}$ is practical.
It is easy to check that
\begin{gather*}
\Phi_{12}(3^{3^k})\ls 2 z_k,\
\Phi_{60}(3^{3^k})\ls 2z_k\Phi_{12}(3^{3^k}),
\\\Phi_{84}(3^{3^k})\ls 2z_k\Phi_{12}(3^{3^k})\Phi_{60}(3^{3^k}),\
\Phi_{420}(3^{3^k})\ls 2 z_k\Phi_{12}(3^{3^k})\Phi_{60}(3^{3^k})\Phi_{84}(3^{3^k}).
\end{gather*}
In light of these and \eqref{2.9}, by applying Lemma \ref{twotimes} four times, we see that
$z_{k+1}$ is practical. This concludes the induction step.

The proof of Theorem 1.4 is now complete. \qed

\end{document}